\documentclass[12pt]{article}
\usepackage{makeidx}
\usepackage{amssymb,amsthm,amsmath,amsfonts,enumerate}
\usepackage{graphicx,caption,subcaption,float,relsize}
\usepackage{url}
\usepackage{hyperref}
\usepackage{rotating}
\usepackage[sort]{natbib}
\usepackage{array}
\usepackage[nodisplayskipstretch]{setspace}
\usepackage{indentfirst}

\setcitestyle{square}
\hypersetup{
  colorlinks = true, urlcolor = cyan, linkcolor = blue, citecolor = red }
\renewenvironment{abstract}
 {\small
  \begin{center}
  \bfseries \abstractname\vspace{-.0em}\vspace{0pt}
  \end{center}
  \list{}{    \setlength{\leftmargin}{0mm}
    \setlength{\rightmargin}{\leftmargin}  }  \item\relax}
 {\endlist}

\allowdisplaybreaks
\newtheorem {theorem}{Theorem}[section]

\newtheorem{corollary}[theorem]{Corollary}

\newtheorem{example}{EXAMPLE}

\newtheorem{lemma}[theorem]{Lemma}

\newtheorem{proposition}[theorem]{Proposition}
\newtheorem{remark}{Remark}[section]

\allowdisplaybreaks

\begin{document}

\title{Symmetrization for high dimensional dependent random variables}
\author{Jonathan B. Hill\thanks{%
Department of Economics, University of North Carolina, Chapel Hill, North
Carolina, \texttt{jbhill@email.unc.edu}; \texttt{https://tarheels.live/jbhill%
}. }\medskip \\
Dept. of Economics, University of North Carolina, Chapel Hill, NC}
\date{{\large This draft:} \today
}
\maketitle

\begin{abstract}
We establish a generic symmetrization property for dependent random
variables $\{x_{t}\}_{t=1}^{n}$ on $\mathbb{R}^{p}$, where $p$ $>>$ $n$ is
allowed. We link $\mathbb{E}\psi (\max_{1\leq i\leq
p}|1/n\sum_{t=1}^{n}(x_{i,t}$ $-$ $\mathbb{E}x_{i,t})|)$ to $\mathbb{E}\psi
(\max_{1\leq i\leq p}|1/n$ $\sum_{t=1}^{n}\eta _{t}(x_{i,t}$ $-$ $\mathbb{E}%
x_{i,t})|)$ for non-decreasing convex $\psi $ $:$ $[0,\infty )$ $\rightarrow 
$ $\mathbb{R}$, where $\{\eta _{t}\}_{t=1}^{n}$ are block-wise independent
random variables, with a remainder term based on high dimensional Gaussian
approximations that need not hold at a high level. Conventional usage of $%
\eta _{t}(x_{i,t}$ $-$ $\tilde{x}_{i,t})$ with $\{\tilde{x}%
_{i,t}\}_{t=1}^{n} $ an independent copy of $\{x_{i,t}\}_{t=1}^{n}$, and
Rademacher $\eta _{t}$, is not required in a generic environment, although
we may trivially replace $\mathbb{E}x_{i,t}$ with $\tilde{x}_{i,t}$. In the
latter case with Rademacher $\eta _{t}$ our result reduces to classic
symmetrization under independence. We bound and therefore verify the
Gaussian approximations in mixing and physical dependence settings, thus
bounding $\mathbb{E}\psi (\max_{1\leq i\leq p}|1/n\sum_{t=1}^{n}(x_{i,t}$ $-$
$\mathbb{E}x_{i,t})|)$; and apply the main result to a generic %
\citet{Nemirovski2000}-like $\mathcal{L}_{q}$-maximal moment bound for $%
\mathbb{E}\max_{1\leq i\leq p}|1/n\sum_{t=1}^{n}(x_{i,t}$ $-$ $\mathbb{E}%
x_{i,t})|^{q}$, $q$ $\geq $ $1$.\medskip \newline
 \smallskip \newline
\end{abstract}

\setstretch{1.75}

\section{Introduction\label{sec:intro}}

Let $\{x_{t}\}_{t=1}^{n}$ be a sample of $\mathbb{R}^{p}$-valued random
variables $x_{t}$ $=$ $[x_{i,t}]_{i=1}^{p}$ on a complete probability space $%
(\Omega ,\mathcal{F},\mathbb{P})$, $p$ $\geq $ $1$, where high
dimensionality $p$ $>>$ $n$ is possible. Let $\psi $ be a non-decreasing
convex function on $[0,\infty )$ with $\psi (0)$ $=$ $0$, and let $\{\eta
_{t}\}_{t=1}^{n}$ be block-wise independent random variables. Write $%
\max_{i} $ $:=$ $\max_{1\leq i\leq p}$. We prove a generic \textquotedblleft
symmetrization\textquotedblright -like result (expectations are assumed to
exist):%
\begin{eqnarray}
\mathbb{E}\psi \left( \max_{i}\left\vert \frac{1}{n}\sum_{t=1}^{n}\left(
x_{i,t}-\mathbb{E}x_{i,t}\right) \right\vert \right) &\leq &\frac{1}{2}%
\mathbb{E}\psi \left( 2\max_{i}\left\vert \frac{1}{n}\sum_{t=1}^{n}\eta
_{t}\left( x_{i,t}-\mathbb{E}x_{i,t}\right) \right\vert \right) +\mathcal{R}%
_{n}(p)  \label{main_convex} \\
&\leq &\mathbb{E}\psi \left( 2\max_{i}\left\vert \frac{1}{n}%
\sum_{t=1}^{n}\left( x_{i,t}-\mathbb{E}x_{i,t}\right) \right\vert \right) +2%
\mathcal{R}_{n}(p).  \notag
\end{eqnarray}%
The appearance of the scales $1/2$ and $2$ in $\frac{1}{2}\mathbb{E}\psi
\left( 2\cdots \right) $, contrary to the classic result, are due to use of
a negligible truncation approximation and convexity of $\psi $.
Symmetrization has been a remarkably powerful tool for bounding norms of
independent random vectors that may otherwise be difficult in the absence of
information. Applications include Donsker theorems, Glivenko-Cantelli
theorems and $\mathcal{L}_{q}$-bounds in high dimension (see, e.g., 
\citet[Chapt.
2.3]{vanderVaartWellner1996} and \cite{Nemirovski2000}).

Our method of proof is completely different than standard symmetrization
arguments under independence (cf. \citet{Pollard1984,vanderVaartWellner1996}%
). We prove (\ref{main_convex}) at a high level under minimal assumptions,
based on high dimensional Gaussian approximation arguments related to the
multiplier (wild) dependent block bootstrap. That said, a Gaussian
approximation for $1/\sqrt{n}\sum_{t=1}^{n}(x_{i,t}$ $-$ $\mathbb{E}x_{i,t})$
is itself not assumed to hold.

Indeed, the remainder $\mathcal{R}_{n}(p)$ is a function of $(i)$ a
diverging truncation point used in an asymptotically negligible truncation
approximation; $(ii)$ Gaussian approximation Kolmogorov distances, with and
without blocking. The latter reduce to an $l_{\infty }$ moment and $\ln
(p)/n $ in a variety of settings, hence in those settings $\mathcal{R}%
_{n}(p) $ $\rightarrow $ $0$ provided $p$ $\rightarrow $ $\infty $ as $n$ $%
\rightarrow $ $\infty $ at a controlled rate that depends on $\psi $ and
tail decay properties. The result carries over to any \textit{Orlicz} norm $%
||X||_{\psi }$ $:=$ $\inf \{c$ $>$ $0$ $:$ $\mathbb{E}\psi (X/c)$ $\leq $ $%
1\}$ by using convexity, non-decreasingness, $||aX||_{\psi }$ $=$ $|a|\times
||X||_{\psi }$ and the triangle inequality. We prove $\lim_{n\rightarrow
\infty }\mathcal{R}_{n}(p)$ $=$ $0$ for mixing and physical dependent random
variables (indeed, $\mathcal{R}_{n}(p)$ $\rightarrow $ $0$ for \textit{any}
dependent random variables for which a \textit{negligible} high dimensional
Gaussian approximation exists). Thus as $n$ $\rightarrow $ $\infty $ we get
the usual symmetrization and desymmetrization,%
\begin{eqnarray*}
\mathbb{E}\psi \left( \max_{i}\left\vert \frac{1}{n}\sum_{t=1}^{n}\left(
x_{i,t}-\mathbb{E}x_{i,t}\right) \right\vert \right) &\leq &\frac{1}{2}%
\mathbb{E}\psi \left( 2\max_{i}\left\vert \frac{1}{n}\sum_{t=1}^{n}\eta
_{t}\left( x_{i,t}-\mathbb{E}x_{i,t}\right) \right\vert \right) \\
&\leq &\mathbb{E}\psi \left( 2\max_{i}\left\vert \frac{1}{n}%
\sum_{t=1}^{n}\left( x_{i,t}-\mathbb{E}x_{i,t}\right) \right\vert \right) .
\end{eqnarray*}

Allowing for arbitrary dependence is a boon for broad applicability. In
social and material sciences dependence structures of observed processes are
generally unknown. Furthermore, high dimensionality is encountered in many
disciplines due to the massive amount of data used, arising from survey
techniques and available technology for data collection. Examples span
social, communication, bio-genetic, electrical, and engineering sciences to
name a few: see, e.g., \citet{FanLi2006}, \citet{BuhlmannVanDeGeer2011}, %
\citet{FanLvQi2011}, and \citet{BelloniChernozhukovHansen2014}.

Use of an independent copy $\tilde{x}_{i,t}$ and Rademacher $\eta _{t}$ 
\textit{per se} under dependence do not expedite the proof as it does in the
classic independence setting (e.g. 
\citet[Chapt.
2.3.1-2.3.2]{vanderVaartWellner1996}). Recall $\eta _{t}$\ is Rademacher
when $\mathbb{P}(\eta _{t}$ $=$ $-1)$ $=$ $\mathbb{P}(\eta _{t}$ $=$ $1)$ $=$
$1/2$. Indeed, at a high level we do not impose any additional structure on $%
\eta _{t}$ other than block-wise independence. We only require properties
when we verify $\mathcal{R}_{n}(p)$ $\rightarrow $ $0$, or apply the results
to a high dimensional $\mathcal{L}_{q}$-moment bound. In such cases we
assume $\eta _{t}$ is bounded, while $\mathcal{R}_{n}(p)$ $\rightarrow $ $0$
requires $\mathbb{E}\eta _{t}^{2}$ $=$ $1$, although at the expense of more
intense notation we could assume a general sub-exponential tail structure.
Thus, the developed results under dependence suggests symmetrization \textit{%
in spirit}. However, it is easily shown that (\ref{main_convex}) holds with $%
\eta _{t}(x_{i,t}$ $-$ $\mathbb{E}x_{i,t})$ replaced with $\eta _{t}(x_{i,t}$
$-$ $\tilde{x}_{i,t})$ for Rademacher $\eta _{t}$, translating to classic
symmetrization with a caveat: under dependence while $x_{i,t}$ $-$ $\tilde{x}%
_{i,t}$ is symmetrically distributed and has the same distribution as $\eta
_{t}(x_{i,t}$ $-$ $\tilde{x}_{i,t})$, $\max_{i}|1/n\sum_{t=1}^{n}(x_{i,t}$ $%
- $ $\tilde{x}_{i,t})|$ and $\max_{i}|1/n\sum_{t=1}^{n}\eta _{t}(x_{i,t}$ $-$
$\tilde{x}_{i,t})|$ generally do \textit{not} have the same distribution.
Thus we cannot conclude equality $\mathbb{E}\psi
(\max_{i}|1/n\sum_{t=1}^{n}(x_{i,t}$ $-$ $\tilde{x}_{i,t})|)$ $=$ $\mathbb{E}%
\psi (\max_{i}|1/n\sum_{t=1}^{n}\eta _{t}(x_{i,t}$ $-$ $\tilde{x}_{i,t})|)$
as we do under independence with iid Rademacher $\eta _{t}$ (%
\citet[p.
109]{vanderVaartWellner1996}). Thus we have a remainder term $\mathcal{R}%
_{n}(p)$. However, $\mathcal{R}_{n}(p)$ $=$ $0$ under independence,
rendering classic symmetrization in that case.

A major purpose of symmetrization is to make it possible to bound norms of
random vectors in the absence of good control on the distribution. But such
an absence is only for independent $x_{i,t}$, eventually under some higher
moment condition (depending on how symmetrization is used, e.g. moment
bound). Thus there is to date always assumed joint distribution control, 
\textit{independence}, which with only mild additional assumptions implies $%
1/\sqrt{n}\sum_{t=1}^{n}(x_{i,t}$ $-$ $\mathbb{E}x_{i,t})$ belongs to the
domain of attraction of a normal law. In this paper we free-up that control
by permitting dependent and heterogeneous data, which necessitates the use
of blocking as discussed above.

In Section \ref{sec:symm} we prove (\ref{main_convex}) by using telescoping
blocks, first assuming $x_{t}$ is bounded. We subsequently allow $x_{t}$ to
be unbounded by using a truncation approximation that builds on results
under boundedness. We apply the main result in Section \ref{sec:max_mom} to
a new maximal moment inequality in the style of \citet{Nemirovski2000},
except instead of independence we allow for physical dependence as in %
\citet{Wu2005} and \citet{WuMin2005}. The appendix contains omitted proofs.
Finally, examples in which the required high dimensional Gaussian
approximations are negligible are presented in the supplemental appendix %
\citet[Appendix B]{symm_sm}.

Throughout $\{x_{t}\}_{t\in \mathbb{N}}$ have non-degenerate distributions. $%
\mathbb{E}$ is the expectations operator; $\mathbb{E}_{\mathcal{A}}$ is
expectations conditional on $\mathcal{F}$-measurable $\mathcal{A}$. $%
\mathcal{L}_{q}$ $:=$ $\{X,$ $\sigma (X)$ $\subset $ $\mathcal{F}:$ $\mathbb{%
E}|X|^{q}$ $<$ $\infty \}$. $||\cdot ||_{q}$ is the $\mathcal{L}_{q}$-norm. $%
a.s.$ is $\mathbb{P}$-\textit{almost surely}. $K$ $>$ $0$ is a finite
constant that may have different values in different places. Similarly
infinitessimal $\iota $ $>$ $0$ may change line to line. $x$ $\lesssim $ $y$
if $x$ $\leq $ $Ky$ for some $K$ $>$ $0$ that is not a function of $n$.
Similarly $x$ $\simeq $ $y$ if $x/y$ $\rightarrow $ $K$ $>$ $0$. Write $%
\max_{i}$ $:=$ $\max_{1\leq i\leq p}$ and $\max_{i,t}$ $:=$ $\max_{1\leq
i\leq p}\max_{1\leq t\leq n}$.

\section{Symmetrization\label{sec:symm}}

Let $\Psi $ be the class of non-decreasing convex functions that are
continuously differentiable on their support: 
\begin{equation*}
\Psi :=\left\{ \psi :[0,\infty )\rightarrow \mathbb{R}:\psi (x)\leq \psi (y)%
\text{ }\forall y\geq x\text{ and }\psi (0)=0\right\} .
\end{equation*}%
Classic examples include the $l_{q}$-metric $x^{q}$\ and the centered
exponential $\exp \{ax^{b}\}$ $-$ $1$, $(a,b)$ $>$ $0$, for $x$ $\geq $ $0$.
Continuous differentiability with nondecreasingness yields a (generalized)
inverse function which we exploit for expectations computation. We can do
away with differentiability by using a well known bound for convex functions 
$\mathbb{E}\psi \left( \left\vert X\right\vert \right) $ $\leq $ $[\psi
\left( b\right) /b]\mathbb{E}|X|$ when $\mathbb{P}(X$ $\in $ $[-b,b])$ $=$ $%
1 $, with $\psi $ $:$ $U$ $\rightarrow $ $\mathbb{R}$, $[0,b]$ $\subseteq $ $%
U$ with $\psi (0)$ $=$ $0$ (see \citet{Edmundson1956} and %
\citet{Madansky1959}). Assume throughout $\mathbb{E}x_{t}$ $=$ $0$.

\subsection{Dependence: bounded\label{sec:dep_b}}

We initially assume $\{x_{t}\}_{t\in \mathbb{N}}$ are bounded variables on $%
\mathbb{R}^{p}$, $p$ $\geq $ $1$. We then use a truncation approximation
that relies on arguments under boundedness.

In order to \textquotedblleft symmetrize\textquotedblright\ with an iid
multiplier that yields the same dependence structure as $\{x_{t}\}_{t=1}^{n}$
asymptotically, we use expanding sub-sample blocks and block-wise
independent multipliers (cf. \citet{Liu1988,Kunsch1989,PolitisRomano1994}).
Let $b_{n}$ $\in $ $\{1,...,n$ $-$ $1\}$ be a pre-set block size, $b_{n}$ $%
\rightarrow $ $\infty $, $b_{n}$ $=$ $o(n)$. Define $\mathcal{N}_{n}$ $:=$ $%
[n/b_{n}]$, and index sets $\mathfrak{B}_{l}$ $:=$ $\{(l$ $-$ $1)b_{n}$ $+$ $%
1,\dots ,lb_{n}\}$ with $l$ $=$ $1,\dots ,\mathcal{N}_{n}$, and assume $%
\mathcal{N}_{n}b_{n}$ $=$ $n$ throughout to reduce notation. Generate
independent random numbers $\{\varepsilon _{l}\}_{l=1}^{\mathcal{N}_{n}}$,
and define the sample $\{\eta _{t}\}_{t=1}^{n}$ by setting $\eta _{t}$ $=$ $%
\varepsilon _{l}$ if $t$ $\in $ $\mathfrak{B}_{l}$. Define%
\begin{eqnarray}
&&\mathcal{X}_{n}(i):=\frac{1}{\sqrt{n}}\sum_{t=1}^{n}x_{i,t}  \label{X_sig}
\\
&&\mathcal{X}_{n}^{\ast }(i):=\frac{1}{\sqrt{n}}\sum_{t=1}^{n}\eta
_{t}x_{i,t}=\frac{1}{\sqrt{n}}\sum_{l=1}^{\mathcal{N}_{n}}\varepsilon _{l}%
\mathcal{S}_{n,l}(i)\text{ where }\mathcal{S}_{n,l}(i):=%
\sum_{t=(l-1)b_{n}+1}^{lb_{n}}x_{i,t}.\text{ \ \ \ \ \ \ \ }  \notag
\end{eqnarray}%
Only \textquotedblleft big\textquotedblright\ blocks $\mathcal{S}_{n,l}(i)$
are used here. In comparable high dimensional settings see, e.g., %
\citet{Chernozhukov_etal2019} who use big and little blocks, and %
\citet{ZhangCheng2018} who use two mutually independent iid multipliers
separately for big and small blocks. See also \citet{Shao2011}. Any such
related approach can be used here.

Let $\{\boldsymbol{X}_{n}(i)\}_{i=1}^{p}$ be a Gaussian process, $%
\boldsymbol{X}_{n}(i)$ $\sim $ $N(0,\mathbb{E}\mathcal{X}_{n}^{2}(i))$, and
define Gaussian approximation Kolmogorov distances with and without blocking%
\begin{eqnarray}
&&\rho _{n}:=\sup_{z\geq 0}\left\vert \mathbb{P}\left( \max_{i}\left\vert 
\mathcal{X}_{n}(i)\right\vert \leq z\right) -\mathbb{P}\left(
\max_{i}\left\vert \boldsymbol{X}_{n}(i)\right\vert \leq z\right) \right\vert
\label{rhos} \\
&&\rho _{n}^{\ast }:=\sup_{z\geq 0}\left\vert \mathbb{P}\left(
\max_{i}\left\vert \mathcal{X}_{n}^{\ast }(i)\right\vert \leq z\right) -%
\mathbb{P}\left( \max_{i}\left\vert \boldsymbol{X}_{n}(i)\right\vert \leq
z\right) \right\vert .  \notag
\end{eqnarray}%
All that follows carries over to the case where an independent copy $\{%
\tilde{x}_{t}\}_{t=1}^{n}$ of $\{x_{t}\}_{t=1}^{n}$ is used. As discussed in
the introduction, however, we generally gain nothing by using a independent
copy under general dependence. Finally, for some sequence of positive real
numbers $\{\mathcal{U}_{n}\}$ to be defined below, define remainder terms
with $\psi ^{\prime }(u)$ $:=$ $(\partial /\partial u)\psi (u)$,%
\begin{equation}
\mathcal{R}_{n}:=\frac{1}{\sqrt{n}}\left\{ \rho _{n}+\rho _{n}^{\ast
}\right\} \times \int_{0}^{\sqrt{n}\mathcal{U}_{n}}\psi ^{\prime }\left( v/%
\sqrt{n}\right) dv.  \label{Rn}
\end{equation}

\begin{remark}
\normalfont At this level of generality we do not impose any structure on $%
\{\varepsilon _{l}\}_{l=1}^{\mathcal{N}_{n}}$ beyond independence, and we do
not impose asymptotic Gaussian approximations \`{a} la $(\rho _{n},\rho
_{n}^{\ast })$ $\rightarrow $ $0$. That said, for a very broad array of
stochastic processes, $|\mathcal{X}_{n}(i)$ $-$ $\boldsymbol{X}_{n}(i)|$ $%
\overset{p}{\rightarrow }$ $0$ and $|\mathcal{X}_{n}^{\ast }(i)$ $-$ $%
\boldsymbol{X}_{n}(i)|$ $\overset{p}{\rightarrow }$ $0$, and indeed $%
\mathcal{R}_{n}$ $\rightarrow $ $0$. Examples are presented in \cite{symm_sm}%
. This is a necessary trade-off: we achieve asymptotic symmetrization for
any dependent and heterogeneous process that satisfies a Gaussian
approximation. Currently, however, symmetrization holds for any independent
random variable (sans Gaussian approximation that typically holds anyway
under mild additional conditions, cf. \cite{Chernozhukov_etal2013}).
\end{remark}

\begin{proposition}[\textquotedblleft Symmetrization\textquotedblright :\
Dependence, Bounded]
\label{prop:symm_dep_convex_bound}Let $\{\mathcal{U}_{n}\}$ be a sequence of
positive real numbers. Let $\{x_{t}\}_{t\in \mathbb{N}}$ be random variables
on $[-\mathcal{U}_{n},\mathcal{U}_{n}]^{p}$, $p$ $\geq $ $1$, and let $%
\{\varepsilon _{l}\}_{l=1}^{\mathcal{N}_{n}}$ be independent random
variables, independent of $\{x_{t}\}_{t=1}^{n}$. We have 
\begin{equation*}
\mathbb{E}\psi \left( \max_{i}\left\vert \bar{x}_{i,n}\right\vert \right)
\leq \mathbb{E}\psi \left( \max_{i}\left\vert \frac{1}{n}\sum_{l=1}^{%
\mathcal{N}_{n}}\varepsilon _{l}\mathcal{S}_{n,l}(i)\right\vert \right) +%
\mathcal{R}_{n}\leq \mathbb{E}\psi \left( \max_{i}\left\vert \bar{x}%
_{i,n}\right\vert \right) +2\mathcal{R}_{n}.
\end{equation*}
\end{proposition}

\begin{remark}
\normalfont$\mathcal{R}_{n}$ captures the error from using a block-wise
multiplier $\varepsilon _{l}$ under general dependence. If $\psi (x)$ $=$ $%
x^{q}$, $x$ $\geq $ $0$ and $q$ $\geq $ $1$, then $\psi ^{\prime }(v/\sqrt{n}%
)$ $=$ $qn^{-(q-1)/2}v^{q-1}$, hence $\mathcal{R}_{n}$ $=$ $\mathcal{U}%
_{n}^{q}n^{-q/2}\{\rho _{n}$ $+$ $\rho _{n}^{\ast }\}$. Thus $\mathcal{R}%
_{n} $ $=$ $o(1/g_{n})$ for some $g_{n}$ $\rightarrow $ $\infty $ as soon as 
$\rho _{n}\vee \rho _{n}^{\ast }$ $=$ $o(n^{q/2}/[\mathcal{U}_{n}^{q}g_{n}])$%
.
\end{remark}

\begin{remark}
\normalfont It is clear from the proof that $\mathbb{E}\psi
(\max_{i}\left\vert \bar{x}_{i,n}\right\vert )$ $\leq $ $\mathbb{E}\psi
(\max_{i}|1/n\sum_{l=1}^{\mathcal{N}_{n}}\varepsilon _{l}\mathcal{S}%
_{n,l}(i)|)$ $+$ $\mathcal{\breve{R}}_{n}$ $\leq $ $\mathbb{E}\psi
(\max_{i}\left\vert \bar{x}_{i,n}\right\vert )$ $+$ $2\mathcal{\breve{R}}%
_{n} $, where%
\begin{eqnarray*}
&&\mathcal{\breve{R}}_{n}:=\frac{1}{\sqrt{n}}\sup_{z\geq 0}\left\vert 
\mathbb{P}\left( \max_{i}\left\vert \frac{1}{\sqrt{n}}\sum_{t=1}^{n}x_{i,t}%
\right\vert \leq z\right) -\mathbb{P}\left( \max_{i}\left\vert \frac{1}{%
\sqrt{n}}\sum_{l=1}^{\mathcal{N}_{n}}\varepsilon _{l}\mathcal{S}%
_{n,l}(i)\right\vert \leq z\right) \right\vert \\
&&\text{ \ \ \ \ \ \ \ \ \ \ \ }\times \int_{0}^{\sqrt{n}\mathcal{U}%
_{n}}\psi ^{\prime }\left( v/\sqrt{n}\right) dv
\end{eqnarray*}%
Under independence set the block size $b_{n}$ $=$ $1$, thus $\eta _{t}$ $=$ $%
\varepsilon _{t}$ and $\mathcal{\breve{R}}_{n}$ $=$ $0$ yielding classic
symmetrization.
\end{remark}

\subsection{Dependence: unbounded\label{sec:dep_ub}}

Now let $x_{t}$ be $\mathbb{R}^{p}$ valued. We use the decomposition $%
\left\vert \bar{x}_{i,n}\right\vert $ $=$ $\left\vert \bar{x}%
_{i,n}\right\vert \mathcal{I}_{\left\vert \bar{x}_{i,n}\right\vert \leq 
\mathcal{U}_{n}}$ $+$ $\left\vert \bar{x}_{i,n}\right\vert \mathcal{I}%
_{\left\vert \bar{x}_{i,n}\right\vert >\mathcal{U}_{n}}$, where $\{\mathcal{U%
}_{n}\}$ is a sequence of positive real numbers, $\mathcal{U}_{n}$ $%
\rightarrow $ $\infty $, that will be implicitly restricted below. By
convexity%
\begin{eqnarray}
\mathbb{E}\psi \left( \max_{i}\left\vert \bar{x}_{i,n}\right\vert \right)
&\leq &\frac{1}{2}\mathbb{E}\psi \left( 2\max_{i}\left\vert \bar{x}%
_{i,n}\right\vert \mathcal{I}_{\left\vert \bar{x}_{i,n}\right\vert \leq 
\mathcal{U}_{n}}\right) +\frac{1}{2}\mathbb{E}\psi \left(
2\max_{i}\left\vert \bar{x}_{i,n}\right\vert \mathcal{I}_{\left\vert \bar{x}%
_{i,n}\right\vert >\mathcal{U}_{n}}\right)  \label{EE} \\
&=&\mathfrak{E}_{n,1}+\mathfrak{E}_{n,2}.  \notag
\end{eqnarray}%
By a change of variables $v$ $=$ $\psi ^{-1}(u)/2$, and the existence of an
inverse function $\psi ^{-1}(\cdot )$ by nondecreasingness and continuity of 
$\psi (\cdot )$,%
\begin{equation*}
\mathfrak{E}_{n,1}=\frac{1}{2}\int_{0}^{\psi (2\mathcal{U}_{n})}\mathbb{P}%
\left( \max_{i}\left\vert \bar{x}_{i,n}\right\vert >\frac{1}{2}\psi
^{-1}(u)\right) du=\int_{0}^{\mathcal{U}_{n}}\psi ^{\prime }(2v)\mathbb{P}%
\left( \max_{i}\left\vert \bar{x}_{i,n}\right\vert >v\right) dv.
\end{equation*}%
Since the latter integral is bounded, by arguments in the proof of
Proposition \ref{prop:symm_dep_convex_bound}%
\begin{eqnarray*}
\mathfrak{E}_{n,1} &\leq &\frac{1}{2}\frac{1}{\sqrt{n}}\int_{0}^{\sqrt{n}%
\mathcal{U}_{n}}\psi ^{\prime }\left( 2v/\sqrt{n}\right) \times \mathbb{P}%
\left( \max_{i}\left\vert \sum_{l=1}^{\mathcal{N}_{n}}\varepsilon _{l}\frac{%
\mathcal{S}_{n,l}(i)}{\sqrt{n}}\right\vert \leq v\right) dv \\
&&\text{ \ \ \ \ \ \ \ }+\frac{1}{2}\left\{ \rho _{n}+\rho _{n}^{\ast
}\right\} \frac{1}{\sqrt{n}}\int_{0}^{\sqrt{n}\mathcal{U}_{n}}\psi ^{\prime
}\left( 2v/\sqrt{n}\right) dv \\
&\leq &\frac{1}{2}\mathbb{E}\psi \left( 2\max_{i}\left\vert \frac{1}{n}%
\sum_{l=1}^{\mathcal{N}_{n}}\varepsilon _{l}\mathcal{S}_{n,l}(i)\right\vert
\right) +\mathcal{R}_{n,1}^{\prime }
\end{eqnarray*}%
with blocking induced remainder 
\begin{equation*}
\mathcal{R}_{n,1}^{\prime }:=\frac{1}{2}\{\rho _{n}+\rho _{n}^{\ast }\}\frac{%
1}{\sqrt{n}}\int_{0}^{\sqrt{n}\mathcal{U}_{n}}\psi ^{\prime }\left( 2v/\sqrt{%
n}\right) dv.
\end{equation*}

Consider the second term $\mathfrak{E}_{n,2}$\ in (\ref{EE}). Use $\psi (0)$ 
$=$ $0$ and convexity to deduce 
\begin{equation}
\mathfrak{E}_{n,2}=\frac{1}{2}\mathbb{E}\psi \left( 2\max_{i}\left\vert \bar{%
x}_{i,n}\right\vert \mathcal{I}_{\left\vert \bar{x}_{i,n}\right\vert >%
\mathcal{U}_{n}}\right) \leq \frac{1}{2}\mathbb{E}\left[ \mathcal{I}%
_{\max_{i}\left\vert \bar{x}_{i,n}\right\vert >\mathcal{U}_{n}}\times \psi
\left( 2\max_{i}\left\vert \bar{x}_{i,n}\right\vert \right) \right] .
\label{En2}
\end{equation}%
Hence by H\"{o}lder and Minkowski inequalities, and convexity and
nondecreasingness,%
\begin{eqnarray}
\mathfrak{E}_{n,2} &\leq &\frac{1}{2}\mathbb{P}\left( \max_{i}\left\vert 
\bar{x}_{i,n}\right\vert \geq \mathcal{U}_{n}\right) ^{(r-1)/r}\times
\left\Vert \psi \left( 2\max_{i}\left\vert \bar{x}_{i,n}\right\vert \right)
\right\Vert _{r}\text{ for }r>1  \notag \\
&\leq &\mathcal{R}_{n,2}^{\prime }:=\frac{1}{2}\mathbb{P}\left(
\max_{i}\left\vert \bar{x}_{i,n}\right\vert \geq \mathcal{U}_{n}\right)
^{(r-1)/r}\times \max_{t}\left\Vert \psi \left( 2\max_{i}\left\vert
x_{i,t}\right\vert \right) \right\Vert _{r},  \label{Rn2}
\end{eqnarray}%
with a truncation induced remainder $\mathcal{R}_{n,2}^{\prime }$. This,
along with a standard desymmetrization argument, proves the main result of
the paper.

\begin{proposition}[\textquotedblleft Symmetrization\textquotedblright :\
Dependence, Unbounded]
\label{prop:symm_dep_convex_unbound}Let $\{x_{t}\}_{t\in \mathbb{N}}$ be
random variables on $\mathbb{R}^{p}$, $p$ $\geq $ $1$, and let $%
\{\varepsilon _{l}\}_{l=1}^{\mathcal{N}_{n}}$ be iid random variables,
independent of $\{x_{t}\}_{t=1}^{n}$. Assume $||\psi (2\max_{i}\left\vert
x_{i,t}\right\vert )||_{r}$ $<$ $\infty $ for each $t$ and some $r$ $>$ $1$.
Then%
\begin{eqnarray*}
\mathbb{E}\psi \left( \max_{i}\left\vert \bar{x}_{i,n}\right\vert \right)
&\leq &\frac{1}{2}\mathbb{E}\psi \left( 2\max_{i}\left\vert \frac{1}{n}%
\sum_{l=1}^{\mathcal{N}_{n}}\varepsilon _{l}\mathcal{S}_{n,l}(i)\right\vert
\right) +\mathcal{R}_{n,1}^{\prime }+\mathcal{R}_{n,2}^{\prime } \\
&\leq &\frac{1}{2}\mathbb{E}\psi \left( 2\max_{i}\left\vert \bar{x}%
_{i,n}\right\vert \right) +2\left\{ \mathcal{R}_{n,1}^{\prime }+\mathcal{R}%
_{n,2}^{\prime }\right\} .
\end{eqnarray*}
\end{proposition}

\begin{remark}
\normalfont\label{rm:exp_xbar}If $\psi (x)$ $=$ $x^{q}$, $x$ $>$ $0$ and $q$ 
$\geq $ $1$, then $||\psi (2\max_{i}\left\vert x_{i,t}\right\vert )||_{r}$ $%
< $ $\infty $ in the second remainder $\mathcal{R}_{n,2}^{\prime }$ \emph{if
and only if} $\mathbb{E}\max_{i}\left\vert x_{i,t}\right\vert ^{qr}$ $<$ $%
\infty $ for some $r$ $>$ $1$. In the exponential case $\psi (x)$ $=$ $\exp
\{ax^{b}\}$ $-$ $1$, $(a,b)$ $>$ $0$, it requires sub-exponential tails $%
\mathbb{E}\exp \{ra\max_{i}\left\vert x_{i,t}\right\vert ^{b}\}$ $<$ $\infty 
$.
\end{remark}

\begin{remark}
\normalfont The combined remainders for some $r$ $>$ $1,$ 
\begin{eqnarray}
\mathcal{R}_{n,1}^{\prime }+\mathcal{R}_{n,2}^{\prime } &=&\{\rho _{n}+\rho
_{n}^{\ast }\}n^{-1/2}\int_{0}^{\sqrt{n}\mathcal{U}_{n}}\psi ^{\prime
}\left( 2v/\sqrt{n}\right) dv  \label{R12} \\
&&+\frac{1}{2}\mathbb{P}\left( \max_{i}\left\vert \bar{x}_{i,n}\right\vert
\geq \mathcal{U}_{n}\right) ^{(r-1)/r}\max_{t}\left\Vert \psi \left(
2\max_{i}\left\vert x_{i,t}\right\vert \right) \right\Vert _{r},  \notag
\end{eqnarray}%
capture approximation errors from blocking and truncation, respectively. $%
\mathcal{R}_{n,1}^{\prime }$ is monotonically increasing as the truncation
level $\mathcal{U}_{n}$ $\rightarrow $ $\infty $, a penalty for having
dependent (hence blocked) data and thus having Gaussian approximations $%
(\rho _{n},\rho _{n}^{\ast })$. The truncation error $\mathcal{R}%
_{n,2}^{\prime }$, however, is logically monotonically decreasing in $%
\mathcal{U}_{n}$.

Consider an $l_{q}$ map $\psi (x)$ $=$ $x^{q}$, and assume sub-exponential
tails for $|\bar{x}_{i,n}|$ 
\begin{equation}
\mathbb{P(}\left\vert \bar{x}_{i,n}\right\vert >x)=a\exp \{-bn^{\gamma
}x^{\gamma }\mathcal{\}}\text{ }\forall x>0,\text{ }a,b,\gamma >0.
\label{sube}
\end{equation}%
Use Lemma \ref{lm:concen}.c below for $\mathbb{P}(\max_{i}|\bar{x}_{i,n}|$ $%
\geq $ $\mathcal{U}_{n})$ to yield for any $\phi $ $\in $ $(0,\gamma )$ 
\begin{equation*}
\mathcal{R}_{n,1}^{\prime }+\mathcal{R}_{n,2}^{\prime }\lesssim
2^{q-1}\left\{ \rho _{n}+\rho _{n}^{\ast }\right\} \mathcal{U}%
_{n}^{q}+2^{q-1}\left( \frac{\ln (p)}{n^{\phi }\mathcal{U}_{n}^{\phi }\ln
(\ln p)}\right) ^{(r-1)/r}\max_{t}\left\Vert \max_{i}\left\vert
x_{i,t}\right\vert \right\Vert _{qr}^{q}.
\end{equation*}%
Now let $\mathcal{U}_{n}^{\ast }$ minimize the upper bound, thus%
\begin{equation*}
\mathcal{U}_{n}^{\ast }=\left\{ \frac{\phi }{q}\left( \frac{r-1}{r}\right)
\left( \frac{1}{\rho _{n}+\rho _{n}^{\ast }}\right) \left( \frac{\ln (p)}{%
n^{\phi }\ln (\ln p)}\right) ^{(r-1)/r}\max_{t}\left\Vert \max_{i}\left\vert
x_{i,t}\right\vert \right\Vert _{qr}^{q}\right\} ^{\frac{1}{q+\phi (r-1)/r}}.
\end{equation*}%
Faster Gaussian approximation convergence $(\rho _{n},\rho _{n}^{\ast })$ $%
\rightarrow $ $0$ implies truncation-based $\mathcal{R}_{n,2}^{\prime }$
dominates, thus a \emph{larger} truncation point $\mathcal{U}_{n}^{\ast }$\
is best. Conversely, larger $\gamma $ implies thinner tails which admit a
larger nuisance term $\phi $, thus $\mathcal{R}_{n,1}^{\prime }$ dominates.
In this case \emph{smaller} $\mathcal{U}_{n}^{\ast }$\ is best.
\end{remark}

\begin{remark}
\normalfont Notice (\ref{sube}) effectively represents a Bernstein or
Fuk-Naegev-type inequality. The condition is valid when $x_{i,t}$ has
sub-exponential tails and, for example, is physical dependent (%
\citet[Theorem 2($ii$)]{Wu2005}), geometric $\tau $-mixing (%
\citet[Theorem
1]{Merlevede_etal2011}), or $\alpha $-mixing or a mixingale (%
\citet{Hill2024_maxlln,Hill_mixg})
\end{remark}

\begin{remark}
\normalfont In \citet[Appendix B]{symm_sm} we prove $(\rho _{n},\rho
_{n}^{\ast })$ $\rightarrow $ $0$ with bounds on $p$ under mixing and
physical dependence, and a variety of tail conditions.
\end{remark}

Remainder $\mathcal{R}_{n,2}^{\prime }$ in (\ref{Rn2}) has a tail measure $%
\mathbb{P}(\max_{i}|\bar{x}_{i,n}|$ $\geq $ $\mathcal{U}_{n})$. Besides
classic concentration bounds like the union bound with Markov's or
Chernoff's inequality, this can be bounded in a variety of ways, akin to
Nemirovski's bound (\citet{Nemirovski2000,BuhlmannVanDeGeer2011}). Define $%
\mathbb{\bar{P}}_{\mathcal{U}}$ $:=$ $\max_{i}\mathbb{P}(|\bar{x}_{i,n}|$ $%
\geq $ $\mathcal{U})$ for any $\mathcal{U}>0$.

\begin{lemma}
\label{lm:concen}Let $\{x_{t}\}$ be random variables on $\mathbb{R}^{p}$.$%
\medskip $\newline
$a.$ In general $\mathbb{P}\left( \max_{i}\left\vert \bar{x}%
_{i,n}\right\vert \geq \mathcal{U}_{n}\right) $ $\lesssim $ $2\ln (p)/\ln (%
\mathbb{\bar{P}}_{\mathcal{U}_{n}}^{-1}\ln (p))$.$\medskip $\newline
$b.$ If $x_{i,t}$ are $\mathcal{L}_{q}$-bounded, $q$ $\geq $ $1$, then $%
\mathbb{P}\left( \max_{i}\left\vert \bar{x}_{i,n}\right\vert \geq \mathcal{U}%
_{n}\right) $ $\lesssim $ $2\ln (p)/\ln (\mathcal{U}_{n}^{q}[\max_{i}\mathbb{%
E}|\bar{x}_{i,n}|^{q}]^{-1}\ln (p)).\medskip $\newline
$c.$ If $\mathbb{P}(|\bar{x}_{i,n}|$ $\geq $ $c)$ $\leq $ $a\exp
\{-bn^{\gamma }c^{\gamma }\}$ $\forall c$ $>$ $0$ and some $a,b,\gamma $ $>$ 
$0$, then $\mathbb{P}\left( \max_{i}\left\vert \bar{x}_{i,n}\right\vert \geq 
\mathcal{U}_{n}\right) \lesssim \ln (p)/[n^{\phi }\mathcal{U}_{n}^{\phi }\ln
(\ln p)]$ for any $\phi $ $\in $ $(0,\gamma )$, $p$ $>$ $e$ and $\ln (p)$ $%
\lesssim $ $\exp \{\mathcal{K}n^{\gamma -\phi }\mathcal{U}_{n}^{\gamma -\phi
}\}$ for all $\mathcal{K}$ $>$ $0.$
\end{lemma}

\begin{remark}
\normalfont We use a conventional log-exp bound with tuning parameter $%
\lambda $ $>$ $0$ in order to prove the claims. ($a$) optimizes the bound
without use of higher moments, while ($b$) and ($c$) optimize the bound with
higher moments,, cf. Remark \ref{rm:exp_xbar}.
\end{remark}

\begin{remark}
\normalfont The condition $\ln (p)$ $\lesssim $ $\exp \{\mathcal{K}n^{\gamma
-\phi }\mathcal{U}_{n}^{\gamma -\phi }\}$ in ($c$) is non-binding
considering $\ln (p)$ $=$ $o(n^{\phi }\mathcal{U}_{n}^{\phi })$ is required
for $\mathbb{P}(\max_{i}|\bar{x}_{i,n}|$ $\geq $ $\mathcal{U}_{n})$ $%
\rightarrow $ $0$.
\end{remark}

\begin{example}[Sub-exponential]
\normalfont Consider $\psi (x)$ $=$ $x^{q}$, $q$ $\geq $ $1$, write $%
\mathcal{M}_{n}$ $:=$ $\max_{t}\left\Vert \max_{i}\left\vert
x_{i,t}\right\vert \right\Vert _{qr}$ and revisit total remainder (\ref{R12}%
) to yield under sub-exponential ($c$) 
\begin{equation*}
\mathcal{R}_{n,1}^{\prime }+\mathcal{R}_{n,2}^{\prime }\lesssim
2^{q-1}\{\rho _{n}+\rho _{n}^{\ast }\}\mathcal{U}_{n}^{q}+\frac{2^{q-1}}{%
\mathcal{U}_{n}^{\phi (r-1)/r}}\left( \frac{\ln (p)}{n^{\phi }\ln (\ln p)}%
\right) ^{(r-1)/r}\mathcal{M}_{n}^{q}.
\end{equation*}%
The upper-bound is minimized with%
\begin{equation*}
\mathcal{U}_{n}^{\ast }=\left\{ \frac{\phi (r-1)}{qr\left\{ \rho _{n}+\rho
_{n}^{\ast }\right\} }\left( \frac{\ln (p)}{n^{\phi }\ln (\ln p)}\right)
^{(r-1)/r}\mathcal{M}_{n}^{q}\right\} ^{\frac{1}{q+\phi (r-1)/r}}.
\end{equation*}%
Thus for some function $\mathcal{K}\left( \phi ,r,q\right) $ $>$ $0$, 
\begin{equation*}
\mathcal{R}_{n,1}^{\prime }+\mathcal{R}_{n,2}^{\prime }\lesssim \mathcal{K}%
\left( \phi ,r,q\right) \{\rho _{n}+\rho _{n}^{\ast }\}^{\frac{\phi (r-1)/r}{%
q+\phi (r-1)/r}}\left\{ \left( \frac{\ln (p)}{n^{\phi }\ln (\ln p)}\right)
^{(r-1)/r}\mathcal{M}_{n}^{q}\right\} ^{\frac{q}{q+\phi (r-1)/r}}.
\end{equation*}%
We naturally need $\{\rho _{n},\rho _{n}^{\ast }\}$ $\rightarrow $ $0$ to
ensure $\mathcal{U}_{n}^{\ast }$ $\rightarrow $ $\infty $ and $\mathcal{R}%
_{n,j}^{\prime }$ $\rightarrow $ $0$. If, for example, $\{\rho _{n},\rho
_{n}^{\ast }\}$ $=$ $o(n^{-\rho })$, $\rho $ $>$ $0$, then $\mathcal{R}%
_{n,1}^{\prime }$ $+$ $\mathcal{R}_{n,2}^{\prime }$ $\rightarrow $ $0$
sufficiently when $\ln (p)$ $=$ $O(n^{\phi (1+\rho /q)}/\mathcal{M}%
_{n}^{qr/(r-1)})$. See \citet[Appendix B]{supp_mat} for conditions yielding $%
\{\rho _{n},\rho _{n}^{\ast }\}$ $=$ $o(n^{-\rho })$.
\end{example}

\section{Application: maximal moment inequality\label{sec:max_mom}}

We apply the main result to deduce a new maximal moment inequality. Set
throughout $\psi (x)$ $=$ $x^{q}$, $x$ $\geq $ $0$ and $q$ $\geq $ $1$. The
following mimics classic arguments based on (conditional) Hoeffding's
inequality, here extended to block-wise partial sums. Let $\{\varepsilon
_{l}\}_{l=1}^{\mathcal{N}_{n}}$ be iid, zero mean and bounded $\mathbb{P}%
(|\varepsilon _{l}|$ $<$ $c)$ $=$ $1$ for some $c$ $\in $ $(0,\infty )$.
Write $\mathfrak{X}^{(n)}$ $:=$ $\{x_{t}\}_{t=1}^{n}$. By Jensen and
Hoeffding inequalities (\citet[Lemma
14.14]{BuhlmannVanDeGeer2011})%
\begin{eqnarray*}
\mathbb{E}\max_{i}\left\vert \frac{1}{n}\sum_{l=1}^{\mathcal{N}%
_{n}}\varepsilon _{l}\mathcal{S}_{n,l}(i)\right\vert ^{q} &=&\mathbb{EE}_{%
\mathfrak{X}^{(n)}}\max_{i}\left\vert \frac{1}{n}\sum_{l=1}^{\mathcal{N}%
_{n}}\varepsilon _{l}\mathcal{S}_{n,l}(i)\right\vert ^{q} \\
&\leq &2^{q/2}c^{q}\left( \frac{\ln \left( 2p\right) }{n}\right) ^{q/2}%
\mathbb{E}\left( \max_{i}\left\vert \frac{1}{n}\sum_{l=1}^{\mathcal{N}_{n}}%
\mathcal{S}_{n,l}^{2}(i)\right\vert ^{q/2}\right) .
\end{eqnarray*}%
Recall $c$ $=$ $1$ under the classic Rademacher assumption. See %
\citet{Bentkus2004,Bentkus2008} for generalizations of Hoeffding's
inequality to unbounded $\{\varepsilon _{l}\}_{l=1}^{\mathcal{N}_{n}}$. In
order eventually to achieve negligible remainders $\mathcal{R}_{n,j}^{\prime
}$ $\rightarrow $ $0$ we require $\mathbb{E}\varepsilon _{l}^{2}$ $=$ $1$
for a Gaussian-to-Gaussian comparison; cf. \citet[Appendix
B]{symm_sm}. Thus, not surprisingly $c$ cannot be arbitrarily small.

The preceding with Proposition \ref{prop:symm_dep_convex_unbound} and Lemma %
\ref{lm:concen} prove the following maximal moment inequality. It is
essentially a generalization of \citet{Nemirovski2000}'s moment bound to
otherwise arbitrary random variables by generating remainder terms based on
blocking and negligible truncation.

\begin{theorem}
\label{th:max_mom}Assume $\mathcal{M}_{n}$ $:=$ $\max_{t}||\max_{i}\left%
\vert x_{i,t}\right\vert ||_{qr}$ $<$ $\infty $ for some $r$ $>$ $1$ and
each $n$, where $\mathcal{M}_{n}$ $\rightarrow $ $\infty $ is possible. Let $%
\{\mathcal{U}_{n}\}$ be a sequence of positive real numbers, $\mathcal{U}%
_{n} $ $\rightarrow $ $\infty $. Then for $q$ $\geq $ $1$%
\begin{equation*}
\mathbb{E}\max_{i}\left\vert \bar{x}_{i,n}\right\vert ^{q}\leq
2^{q/2}c^{q}\left( \frac{\ln \left( 2p\right) }{n}\right) ^{q/2}\mathbb{E}%
\left( \max_{i}\left\vert \frac{1}{n}\sum_{l=1}^{\mathcal{N}_{n}}\mathcal{S}%
_{n,l}^{2}(i)\right\vert ^{q/2}\right) +\frac{1}{2}\left\{ \mathcal{R}%
_{n,1}^{\prime }+\mathcal{R}_{n,2}^{\prime }\right\}
\end{equation*}%
where $\mathcal{R}_{n,1}^{\prime }$ $=$ $2^{q}\mathcal{U}_{n}^{q}n^{-q/2}\{%
\rho _{n}+\rho _{n}^{\ast }\}$, and $\mathcal{R}_{n,2}^{\prime }$ is derived
by case as follows.$\medskip $\newline
$a.$ Under $\mathcal{L}_{qr}$-boundedness%
\begin{equation*}
\mathcal{R}_{n,2}^{\prime }\leq 2^{q-1/r}\left( \frac{\ln (p)}{\ln (\mathcal{%
U}_{n}^{q}[\max_{i}\mathbb{E}|\bar{x}_{i,n}|^{q}]^{-1}\ln (p))}\right)
^{(r-1)/r}\mathcal{M}_{n}^{q}.
\end{equation*}%
$b.$ If $\mathbb{P}(|\bar{x}_{i,n}|$ $\geq $ $c)$ $\leq $ $a\exp
\{-bn^{\gamma }c^{\gamma }\}$ $\forall c$ $>$ $0$ for some $a,b,\gamma $ $>$ 
$0$, then for any $\phi $ $\in $ $(0,\gamma )$, $p$ $>$ $e$ and $\ln (p)$ $%
\lesssim $ $\exp \{\mathcal{K}n^{\gamma -\phi }\mathcal{U}_{n}^{\gamma -\phi
}\}$ for all $\mathcal{K}$ $>$ $0$,%
\begin{equation*}
\mathcal{R}_{n,2}^{\prime }\lesssim 2^{q-1}\left( \frac{\ln (p)}{n^{\phi }%
\mathcal{U}_{n}^{\phi }\ln (\ln p)}\right) ^{(r-1)/r}\mathcal{M}_{n}^{q}.
\end{equation*}
\end{theorem}

Theorem \ref{th:max_mom} instantly yields the following.

\begin{corollary}
\label{cor:max_mom}Let the truncation points satisfy $\mathcal{U}_{n}$ $=$ $%
o(n^{1/2})$. Under either of the following settings, for some positive
sequence $\{g_{n}\}$, $g_{n}$ $\rightarrow $ $\infty $, to be implicitly
defined below, and $q$ $\geq $ $1$ 
\begin{equation}
\mathbb{E}\max_{i}\left\vert \bar{x}_{i,n}\right\vert ^{q}\leq
2^{q/2}c^{q}\left( \frac{\ln \left( 2p\right) }{n}\right) ^{q/2}\mathbb{E}%
\left( \max_{i}\frac{1}{n}\sum_{l=1}^{\mathcal{N}_{n}}\mathcal{S}%
_{n,l}^{2}(i)\right) ^{q/2}+o\left( \rho _{n}+\rho _{n}^{\ast }\right)
+o\left( 1/g_{n}\right) .  \label{Eb}
\end{equation}%
$a.$ $x_{i,t}$ are $\mathcal{L}_{qr}$-bounded, $r$ $>$ $1$, and $\ln (p)$ $=$
$o(g_{n}^{-1}\mathcal{M}_{n}^{-qr/(r-1)}\ln [n/\max_{i}\mathbb{E}\left\vert 
\bar{x}_{i,n}\right\vert ^{q}])$.$\medskip $\newline
$b.$ $\mathbb{P}(|\bar{x}_{i,n}|$ $\geq $ $c)$ $\leq $ $a\exp \{-bn^{\gamma
}c^{\gamma }\}$ $\forall c$ $>$ $0$ and some $a,b,\gamma $ $>$ $0$; and for
any $\phi $ $\in $ $(0,\gamma )$, $p$ $>$ $e$, and some $r$ $>$ $1$, we have 
$\ln (p)$ $=$ $o(g_{n}^{-r/(r-1)}\mathcal{M}_{n}^{-qr/(r-1)}n^{3\phi /2})$.
\end{corollary}

\begin{remark}
\normalfont Consider ($a$) and let $\{x_{i,t}\}$ be stationary and uniformly 
$\mathcal{L}_{rq}$-bounded over $i$. A wide array of weak dependence
properties support $\mathbb{E}\left\vert \bar{x}_{i,n}\right\vert ^{q}$ $=$ $%
O(1/n^{q/2})$, including various mixing, mixingale, and physical dependence
(e.g. \citet{Hansen1991,Hansen1992,Wu2005}). Now use $\mathcal{M}_{n}$ $\leq 
$ $p^{1/(rq)}(\max_{i}\mathbb{E}\left\vert x_{i,t}\right\vert ^{rq})^{1/(rq)}
$ to yield $p$ $=$ $o(\{g_{n}^{-1}\ln (n)\}^{r-1})$, thus $g_{n}$ $=$ $o(\ln
(n))$. If cross-coordinate $i$ dependence is known than a potentially vastly
sharper bound on $\mathcal{M}_{n}$ is available. For example, if $\{x_{i,t},%
\mathfrak{F}_{n,i}\}_{i=1}^{k_{n}}$ forms a martingale for some filtration $%
\mathfrak{F}_{n,i}$ then $\mathcal{M}_{n}$ $=$ $O(1)$ for any $p$ by Doob's
inequality. See \cite{Hill2024_maxlln} for examples and theory.
\end{remark}

\begin{remark}
\normalfont Under ($b$) suppose also $\mathbb{P}(|x_{i,t}|$ $\geq $ $c)$ $%
\leq $ $a\exp \{-bc^{\gamma }\}$ for some $\gamma $ $\geq $ $1$, thus $%
\mathcal{M}_{n}$ $=$ $O(\ln (p)^{\psi })$ for some $\psi $ that depends on $%
\gamma ,q,r$. Cf. Remark \ref{rm:exp_xbar}. Moreover, $r$ may be arbitrarily
large under sub-exponential tails, so take $r$ $\rightarrow $ $\infty $.
Therefore $\ln (p)$ $=$ $o(\{n^{3\phi /2}/g_{n}\}^{1/(1+q\psi )})$ and thus $%
g_{n}$ $=$ $o(n^{3\phi /2})$. Now suppose $\gamma $ $=$ $1$ yielding classic
sub-exponential decay. Set $g_{n}$ $=$ $n^{3\phi /4}$ and $\phi $ $=$ $%
\gamma $ $-$ $\iota $ $=$ $1$ $-$ $\iota $ for infinitessimal $\iota $ $>$ $%
0 $\ to yield in (\ref{Eb}) an upper bound remainder $o(\rho _{n}$ $+$ $\rho
_{n}^{\ast }$ $+$ $n^{-3/4+\iota })$ when $\ln (p)$ $=$ $o(n^{3/[4(1+q\psi
)]-\iota })$.
\end{remark}

\subsection{Conclusion}

We extend the symmetrization concept to arbitrarily dependent random
variables by using a negligible truncation approximation, telescoping blocks
with a block-wise dependent multiplier in order to imitate the underlying
dependence structure, and high dimensional Gaussian comparisons. We
therefore sidestep classic arguments utilizing an iid Rademacher multiplier
and independent copy: the multiplier cannot be independent, while the
Rademacher structure serves a far more narrow purpose here (boundedness);
and an independent copy is essentially superfluous under dependence. The
main bound involves remainder terms, errors generated from blocking and
truncation. The multiplier need not be specified at a high level, but will
logically be bounded (or sub-exponential) in applications. We apply the main
result to a new Nemirovski-like moment bound under dependence, and present
examples establishing vanishing Gaussian approximations for mixing and
physical dependent sequences. Future work may focus on sharpness, or utilize
cross-coordinate dependence, issues ignored here for the sake of focus.

\setcounter{equation}{0} \renewcommand{\theequation}{{\thesection}.%
\arabic{equation}} \appendix

\section{Appendix: omitted proofs\label{app:proofs}}

\noindent \noindent \textbf{Proof of Proposition \ref%
{prop:symm_dep_convex_bound}.} The triangle inequality and $\{\rho _{n},\rho
_{n}^{\ast }\}$ defined in (\ref{rhos}) yield%
\begin{equation*}
\sup_{z\geq 0}\left\vert \mathbb{P}\left( \max_{i}\left\vert \frac{1}{\sqrt{n%
}}\sum_{t=1}^{n}x_{i,t}\right\vert \leq z\right) -\mathbb{P}\left(
\max_{i}\left\vert \frac{1}{\sqrt{n}}\sum_{l=1}^{\mathcal{N}_{n}}\varepsilon
_{l}\mathcal{S}_{n,l}(i)\right\vert \leq z\right) \right\vert \leq \rho
_{n}+\rho _{n}^{\ast }.
\end{equation*}%
Replace $x_{i,t}$ with $x_{i,t}/\sqrt{n}$: for each $x$ $\geq $ $0$, 
\begin{eqnarray*}
&&\left\vert \mathbb{P}\left( \max_{i}\left\vert \bar{x}_{i,n}\right\vert
\leq x\right) -\mathbb{P}\left( \max_{i}\left\vert \frac{1}{n}\sum_{l=1}^{%
\mathcal{N}_{n}}\varepsilon _{l}\mathcal{S}_{n,l}(i)\right\vert \leq
x\right) \right\vert \\
&&\text{ \ \ \ \ \ \ \ }=\left\vert \mathbb{P}\left( \max_{i}\left\vert
\sum_{t=1}^{n}\frac{x_{i,t}}{\sqrt{n}}\right\vert \leq \sqrt{n}x\right) -%
\mathbb{P}\left( \max_{i}\left\vert \sum_{l=1}^{\mathcal{N}_{n}}\varepsilon
_{l}\frac{\mathcal{S}_{n,l}(i)}{\sqrt{n}}\right\vert \leq \sqrt{n}x\right)
\right\vert \leq \rho _{n}+\rho _{n}^{\ast }\text{.}
\end{eqnarray*}

Now use a change of variables, the fact that a (generalized) inverse
function $\psi ^{-1}(\cdot )$ exists by nondecreasingness and continuity of $%
\psi (\cdot )$, and $\max_{i}\mathbb{P}(|x_{i,t}|$ $>$ $\mathcal{U}_{n})$ $=$
$0$ to yield for any $q$ $\geq $ $1$%
\begin{eqnarray}
\mathbb{E}\psi \left( \max_{i}\left\vert \bar{x}_{i,n}\right\vert \right)
&=&\int_{0}^{\psi (\mathcal{U}_{n})}\mathbb{P}\left( \psi \left(
\max_{i}\left\vert \bar{x}_{i,n}\right\vert \right) >u\right) du  \notag \\
&=&\int_{0}^{\psi (\mathcal{U}_{n})}\mathbb{P}\left( \max_{i}\left\vert 
\sqrt{n}\bar{x}_{i,n}\right\vert >\sqrt{n}\psi ^{-1}(u)\right) du  \notag \\
&=&\frac{1}{\sqrt{n}}\int_{0}^{\sqrt{n}\mathcal{U}_{n}}\psi ^{\prime }\left(
v/\sqrt{n}\right) \times \mathbb{P}\left( \max_{i}\left\vert \sqrt{n}\bar{x}%
_{i,n}\right\vert >v\right) dv  \notag \\
&\leq &\frac{1}{\sqrt{n}}\int_{0}^{\sqrt{n}\mathcal{U}_{n}}\psi ^{\prime
}\left( v/\sqrt{n}\right) \times \mathbb{P}\left( \max_{i}\left\vert
\sum_{l=1}^{\mathcal{N}_{n}}\varepsilon _{l}\frac{\mathcal{S}_{n,l}(i)}{%
\sqrt{n}}\right\vert \leq v\right) dv  \notag \\
&&\text{ \ \ \ \ \ \ \ }+\left\{ \rho _{n}+\rho _{n}^{\ast }\right\} \frac{1%
}{\sqrt{n}}\int_{0}^{\sqrt{n}\mathcal{U}_{n}}\psi ^{\prime }\left( v/\sqrt{n}%
\right) dv  \notag \\
&=&\mathbb{E}\psi \left( \max_{i}\left\vert \frac{1}{n}\sum_{l=1}^{\mathcal{N%
}_{n}}\varepsilon _{l}\mathcal{S}_{n,l}(i)\right\vert \right) +\mathcal{R}%
_{n},  \label{EEpsi}
\end{eqnarray}%
where $\mathcal{R}_{n}$ $:=$ $\{\rho _{n}$ $+$ $\rho _{n}^{\ast
}\}n^{-1/2}\int_{0}^{\sqrt{n}\mathcal{U}_{n}}\psi ^{\prime }(v/\sqrt{n})dv$.
The last line follows by reversing the change of variables. Repeat the
argument in reverse to yield similarly%
\begin{eqnarray}
\mathbb{E}\psi \left( \max_{i}\left\vert \frac{1}{n}\sum_{l=1}^{\mathcal{N}%
_{n}}\varepsilon _{l}\mathcal{S}_{n,l}(i)\right\vert \right) &\leq &\frac{1}{%
\sqrt{n}}\int_{0}^{\sqrt{n}\mathcal{U}_{n}}\psi ^{\prime }\left( v/\sqrt{n}%
\right) \mathbb{P}\left( \max_{i}\left\vert \sqrt{n}\bar{x}_{i,n}\right\vert
\leq v\right) dv+\mathcal{R}_{n}  \notag \\
&=&\mathbb{E}\psi \left( \max_{i}\left\vert \bar{x}_{i,n}\right\vert \right)
+\mathcal{R}_{n}.  \label{EEpsi2}
\end{eqnarray}%
Combine (\ref{EEpsi}) and (\ref{EEpsi2}) to conclude as claimed%
\begin{equation*}
\mathbb{E}\psi \left( \max_{i}\left\vert \bar{x}_{i,n}\right\vert \right)
\leq \mathbb{E}\psi \left( \max_{i}\left\vert \frac{1}{n}\sum_{l=1}^{%
\mathcal{N}_{n}}\varepsilon _{l}\mathcal{S}_{n,l}(i)\right\vert \right) +%
\mathcal{R}_{n}\leq \mathbb{E}\psi \left( \max_{i}\left\vert \bar{x}%
_{i,n}\right\vert \right) +2\mathcal{R}_{n}.
\end{equation*}%
$\mathcal{QED}$.\medskip \newline
\textbf{Proof of Lemma \ref{lm:concen}. }Set $\mathbb{P}_{\mathcal{U}_{n}}$ $%
:=$ $\max_{i}\mathbb{P}(|\bar{x}_{i,n}|$ $<$ $\mathcal{U}_{n})$ and $\mathbb{%
\bar{P}}_{\mathcal{U}_{n}}$ $:=$ $\max_{i}\mathbb{P}(|\bar{x}_{i,n}|$ $\geq $
$\mathcal{U}_{n})$. Use Jensen's inequality to deduce for any $\lambda $ $>$ 
$0$ 
\begin{equation}
\mathbb{P}\left( \max_{i}\left\vert \bar{x}_{i,n}\right\vert \geq \mathcal{U}%
_{n}\right) \leq \frac{1}{\lambda }\ln \left( \mathbb{E}\left[ \exp \left\{
\lambda \mathcal{I}_{\max_{i}|\bar{x}_{i,n}|\geq \mathcal{U}_{n}}\right\} %
\right] \right) \leq \frac{1}{\lambda }\ln \left( p\max_{i}\mathbb{E}\left[
\exp \left\{ \lambda \mathcal{I}_{|\bar{x}_{i,n}|\geq \mathcal{U}%
_{n}}\right\} \right] \right) .  \notag
\end{equation}%
By construction $\max_{i}\mathbb{E}[\exp \{\lambda \mathcal{I}_{|\bar{x}%
_{i,n}|\geq \mathcal{U}_{n}}\}]$ $\leq $ $\exp \{\lambda \}\mathbb{\bar{P}}_{%
\mathcal{U}_{n}}$ $+$ $\mathbb{P}_{\mathcal{U}_{n}}$ $\leq $ $\exp \{\lambda
\}\mathbb{\bar{P}}_{\mathcal{U}_{n}}$ $+$ $1$. Now use $\ln (1$ $+$ $x)$ $%
\leq $ $x$ $\forall x$ $\geq $ $0$ to yield 
\begin{equation}
\mathbb{P}\left( \max_{i}\left\vert \bar{x}_{i,n}\right\vert \geq \mathcal{U}%
_{n}\right) \leq \frac{1}{\lambda }\ln (p)+\frac{1}{\lambda }\exp \left\{
\lambda \right\} \times \mathbb{\bar{P}}_{\mathcal{U}_{n}}.  \label{PP}
\end{equation}%
\textbf{Claim (a). }Minimize (\ref{PP}) with respect to $\lambda $ to yield $%
\lambda $ $=$ $\ln (\mathbb{\bar{P}}_{\mathcal{U}_{n}}^{-1}\ln (p))$ as $n$ $%
\rightarrow $ $\infty $. Hence%
\begin{equation*}
\mathbb{P}\left( \max_{i}\left\vert \bar{x}_{i,n}\right\vert \geq \mathcal{U}%
_{n}\right) \lesssim 2\ln (p)/\left[ \ln (\mathbb{\bar{P}}_{\mathcal{U}%
_{n}}^{-1}\ln (p))\right] .
\end{equation*}%
\textbf{Claim (b).} Use Markov's inequality $\mathbb{\bar{P}}_{\mathcal{U}%
_{n}}$ $\leq $ $\mathcal{U}_{n}^{-q}\max_{i}\mathbb{E}|\bar{x}_{i,n}|^{q}$
in (\ref{PP}), and the argument under ($a$) to yield the result.\medskip 
\newline
\textbf{Claim (c). }Let $\mathbb{\bar{P}}_{\mathcal{U}_{n}}$ $\leq $ $a\exp
\{-bn^{\gamma }\mathcal{U}_{n}^{\gamma }\}$. By (\ref{PP}) with $\lambda $ $%
= $ $n^{\phi }\mathcal{U}_{n}^{\phi }\ln (\ln p)$ for any $\phi $ $\in $ $%
(0,\gamma )$ we have under $\ln (p)$ $\lesssim $ $\exp \{\mathcal{K}%
n^{\gamma -\phi }\mathcal{U}_{n}^{\gamma -\phi }\}$ for all $\mathcal{K}$ $>$
$0$, 
\begin{eqnarray*}
\mathbb{P}\left( \max_{i}\left\vert \bar{x}_{i,n}\right\vert \geq \mathcal{U}%
_{n}\right) &\leq &\frac{1}{\lambda }\ln (p)+\frac{1}{\lambda }\exp \left\{
\lambda \right\} a\exp \left\{ -bn^{\gamma }\mathcal{U}_{n}^{\gamma }\right\}
\\
&=&\frac{\ln (p)}{n^{\phi }\mathcal{U}_{n}^{\phi }\ln (\ln p)}+a\frac{\left(
\ln (p)\right) ^{n^{\phi }\mathcal{U}_{n}^{\phi }}}{n^{\phi }\mathcal{U}%
_{n}^{\phi }\exp \left\{ bn^{\gamma }\mathcal{U}_{n}^{\gamma }\right\} \ln
(\ln p)}\lesssim \frac{\ln (p)}{n^{\phi }\mathcal{U}_{n}^{\phi }\ln (\ln p)}.
\end{eqnarray*}%
$\mathcal{QED}$.\newline
\setstretch{.75} 
\bibliographystyle{cas-model2-names.bst}
\bibliography{refs_symm}

\begin{thebibliography}{27}
\expandafter\ifx\csname natexlab\endcsname\relax\def\natexlab#1{#1}\fi
\providecommand{\url}[1]{\texttt{#1}}
\providecommand{\href}[2]{#2}
\providecommand{\path}[1]{#1}
\providecommand{\DOIprefix}{doi:}
\providecommand{\ArXivprefix}{arXiv:}
\providecommand{\URLprefix}{URL: }
\providecommand{\Pubmedprefix}{pmid:}
\providecommand{\doi}[1]{\href{http://dx.doi.org/#1}{\path{#1}}}
\providecommand{\Pubmed}[1]{\href{pmid:#1}{\path{#1}}}
\providecommand{\bibinfo}[2]{#2}
\ifx\xfnm\relax \def\xfnm[#1]{\unskip,\space#1}\fi
\bibitem[{Belloni et~al.(2014)Belloni, Chernozhukov and
  Hansen}]{BelloniChernozhukovHansen2014}
\bibinfo{author}{Belloni, A.}, \bibinfo{author}{Chernozhukov, V.},
  \bibinfo{author}{Hansen, C.}, \bibinfo{year}{2014}.
\newblock \bibinfo{title}{High-dimensional methods and inference on structural
  and treatment effects}.
\newblock \bibinfo{journal}{J. Econom. Perspect.} \bibinfo{volume}{28},
  \bibinfo{pages}{29--50}.
\bibitem[{Bentkus(2004)}]{Bentkus2004}
\bibinfo{author}{Bentkus, V.}, \bibinfo{year}{2004}.
\newblock \bibinfo{title}{On hoeffding's inequalities}.
\newblock \bibinfo{journal}{Ann. Probab.} \bibinfo{volume}{32},
  \bibinfo{pages}{1650--1673}.
\bibitem[{Bentkus(2008)}]{Bentkus2008}
\bibinfo{author}{Bentkus, V.}, \bibinfo{year}{2008}.
\newblock \bibinfo{title}{An extension of the hoeffding inequality to unbounded
  random variables}.
\newblock \bibinfo{journal}{Lith. Math. J.} \bibinfo{volume}{48},
  \bibinfo{pages}{137--157}.
\bibitem[{Buhlmann and van~de Geer(2011)}]{BuhlmannVanDeGeer2011}
\bibinfo{author}{Buhlmann, P.}, \bibinfo{author}{van~de Geer, S.},
  \bibinfo{year}{2011}.
\newblock \bibinfo{title}{Statistics for High-Dimensional Data}.
\newblock \bibinfo{publisher}{Springer}, \bibinfo{address}{Berlin}.
\bibitem[{Chernozhukov et~al.(2013)Chernozhukov, Chetverikov and
  Kato}]{Chernozhukov_etal2013}
\bibinfo{author}{Chernozhukov, V.}, \bibinfo{author}{Chetverikov, D.},
  \bibinfo{author}{Kato, K.}, \bibinfo{year}{2013}.
\newblock \bibinfo{title}{Gaussian approximations and multiplier bootstrap for
  maxima of sums of high-dimensional random vectors}.
\newblock \bibinfo{journal}{Ann. Statist.} \bibinfo{volume}{41},
  \bibinfo{pages}{2786--2819}.
\bibitem[{Chernozhukov et~al.(2019)Chernozhukov, Chetverikov and
  Kato}]{Chernozhukov_etal2019}
\bibinfo{author}{Chernozhukov, V.}, \bibinfo{author}{Chetverikov, D.},
  \bibinfo{author}{Kato, K.}, \bibinfo{year}{2019}.
\newblock \bibinfo{title}{Inference on causal and structural parameters using
  many moment inequalities}.
\newblock \bibinfo{journal}{Rev. Econ. Stud.} \bibinfo{volume}{86},
  \bibinfo{pages}{1867--1900}.
\bibitem[{Edmundson(1956)}]{Edmundson1956}
\bibinfo{author}{Edmundson, H.P.}, \bibinfo{year}{1956}.
\newblock \bibinfo{title}{Bounds on the Expectation of a Convex Functions}.
\newblock \bibinfo{type}{Technical Report} \bibinfo{number}{982}. Rand Corp..
  \bibinfo{address}{Santa Monica}.
\bibitem[{Fan and Li(2006)}]{FanLi2006}
\bibinfo{author}{Fan, J.}, \bibinfo{author}{Li, R.}, \bibinfo{year}{2006}.
\newblock \bibinfo{title}{Statistical challenges with high dimensionality:
  Feature selection in knowledge discovery}, in: \bibinfo{editor}{Sanz-Sole,
  M.}, \bibinfo{editor}{Soria, J.}, \bibinfo{editor}{Varona, J.L.},
  \bibinfo{editor}{Verdera, J.} (Eds.), \bibinfo{booktitle}{Proceedings of the
  International Congress of Mathematicians}, \bibinfo{organization}{European
  Mathematical Society}, \bibinfo{address}{Zurich}. pp.
  \bibinfo{pages}{595--622}.
\bibitem[{Fan et~al.(2011)Fan, Lv and Qi}]{FanLvQi2011}
\bibinfo{author}{Fan, J.}, \bibinfo{author}{Lv, J.}, \bibinfo{author}{Qi, .L.},
  \bibinfo{year}{2011}.
\newblock \bibinfo{title}{Sparse high-dimensional models in economics}.
\newblock \bibinfo{journal}{Annu. Rev. Economics} \bibinfo{volume}{3},
  \bibinfo{pages}{291--317}.
\bibitem[{Hansen(1991)}]{Hansen1991}
\bibinfo{author}{Hansen, B.E.}, \bibinfo{year}{1991}.
\newblock \bibinfo{title}{Strong laws for dependent heterogeneous processes}.
\newblock \bibinfo{journal}{Econometric Theory} \bibinfo{volume}{7},
  \bibinfo{pages}{213--221}.
\bibitem[{Hansen(1992)}]{Hansen1992}
\bibinfo{author}{Hansen, B.E.}, \bibinfo{year}{1992}.
\newblock \bibinfo{title}{Erratum: Strong laws for dependent heterogeneous
  processes}.
\newblock \bibinfo{journal}{Econometric Theory} \bibinfo{volume}{8},
  \bibinfo{pages}{421--422}.
\bibitem[{Hill(2024a)}]{Hill2024_maxlln}
\bibinfo{author}{Hill, J.B.}, \bibinfo{year}{2024}a.
\newblock \bibinfo{title}{Max-laws of large numbers for weakly dependent high
  dimensional arrays with applications}.
\newblock \bibinfo{note}{Technical Report, Dept. of Economics, University of
  North Carolina}.
\bibitem[{Hill(2024b)}]{supp_mat}
\bibinfo{author}{Hill, J.B.}, \bibinfo{year}{2024}b.
\newblock \bibinfo{title}{Supplemental material for ``symmetrization for high
  dimensional dependent random variables''}.
\newblock \bibinfo{note}{Dept. of Economics, UNC}.
\bibitem[{Hill(2025a)}]{Hill_mixg}
\bibinfo{author}{Hill, J.B.}, \bibinfo{year}{2025}a.
\newblock \bibinfo{title}{Mixingale and physical dependence equality with
  applications}.
\newblock \bibinfo{journal}{Stat. Probab. Let.} \bibinfo{volume}{in press}.
\bibitem[{Hill(2025b)}]{symm_sm}
\bibinfo{author}{Hill, J.B.}, \bibinfo{year}{2025}b.
\newblock \bibinfo{title}{Supplemental material for ``symmetrization for high
  dimensional dependent random variables''}.
\newblock \bibinfo{note}{Dept. of Economics, University of North Carolina -
  Chapel Hill}.
\bibitem[{K{\"u}nsch(1989)}]{Kunsch1989}
\bibinfo{author}{K{\"u}nsch, H.R.}, \bibinfo{year}{1989}.
\newblock \bibinfo{title}{The jackknife and the bootstrap for general
  stationary observations}.
\newblock \bibinfo{journal}{Ann. Statist.} \bibinfo{volume}{17},
  \bibinfo{pages}{1217--1241}.
\bibitem[{Liu(1988)}]{Liu1988}
\bibinfo{author}{Liu, R.Y.}, \bibinfo{year}{1988}.
\newblock \bibinfo{title}{Bootstrap procedures under some non-i.i.d. models}.
\newblock \bibinfo{journal}{Ann. Statist.} \bibinfo{volume}{16},
  \bibinfo{pages}{1696--1708}.
\bibitem[{Madansky(1959)}]{Madansky1959}
\bibinfo{author}{Madansky, A.}, \bibinfo{year}{1959}.
\newblock \bibinfo{title}{Bounds on the expectation of a convex function of a
  multivariate random variable}.
\newblock \bibinfo{journal}{Ann. Math. Statist.} \bibinfo{volume}{30},
  \bibinfo{pages}{743--746}.
\bibitem[{Merlevede et~al.(2011)Merlevede, Peligrad and
  Rio}]{Merlevede_etal2011}
\bibinfo{author}{Merlevede, F.}, \bibinfo{author}{Peligrad, M.},
  \bibinfo{author}{Rio, E.}, \bibinfo{year}{2011}.
\newblock \bibinfo{title}{Bernstein inequality and moderate deviations for
  weakly dependent sequences}.
\newblock \bibinfo{journal}{Probab. Theory Rel.} \bibinfo{volume}{151},
  \bibinfo{pages}{435--474}.
\newblock \bibinfo{note}{Volume 5}.
\bibitem[{Nemirovski(2000)}]{Nemirovski2000}
\bibinfo{author}{Nemirovski, A.S.}, \bibinfo{year}{2000}.
\newblock \bibinfo{title}{Topics in nonparametric statistics}, in:
  \bibinfo{editor}{Emery, M.}, \bibinfo{editor}{Nemirovski, A.},
  \bibinfo{editor}{Voiculescu, D.}, \bibinfo{editor}{Bernard, P.} (Eds.),
  \bibinfo{booktitle}{Lectures on Probability Theory and Statistics: Ecole
  d'Ete de Probabilites de Saint-Flour XXVIII - 1998}.
  \bibinfo{publisher}{Springer}, \bibinfo{address}{New York}. volume
  \bibinfo{volume}{1738}, pp. \bibinfo{pages}{87--285}.
\bibitem[{Politis and Romano(1994)}]{PolitisRomano1994}
\bibinfo{author}{Politis, D.N.}, \bibinfo{author}{Romano, J.P.},
  \bibinfo{year}{1994}.
\newblock \bibinfo{title}{The stationary bootstrap}.
\newblock \bibinfo{journal}{J. Amer. Statis. Assoc.} \bibinfo{volume}{89},
  \bibinfo{pages}{1303--1313}.
\bibitem[{Pollard(1984)}]{Pollard1984}
\bibinfo{author}{Pollard, D.}, \bibinfo{year}{1984}.
\newblock \bibinfo{title}{Convergence of Stochastic Processes}.
\newblock \bibinfo{publisher}{Springer Verlag}, \bibinfo{address}{New York}.
\bibitem[{Shao(2011)}]{Shao2011}
\bibinfo{author}{Shao, X.}, \bibinfo{year}{2011}.
\newblock \bibinfo{title}{A bootstrap-assisted spectral test of white noise
  under unknown dependence}.
\newblock \bibinfo{journal}{Journal of Econometrics} \bibinfo{volume}{162},
  \bibinfo{pages}{213--224}.
\bibitem[{van~der Vaart and Wellner(1996)}]{vanderVaartWellner1996}
\bibinfo{author}{van~der Vaart, A.}, \bibinfo{author}{Wellner, J.},
  \bibinfo{year}{1996}.
\newblock \bibinfo{title}{Weak Convergence and Empirical Processes}.
\newblock \bibinfo{publisher}{Springer}, \bibinfo{address}{New York}.
\bibitem[{Wu(2005)}]{Wu2005}
\bibinfo{author}{Wu, W.B.}, \bibinfo{year}{2005}.
\newblock \bibinfo{title}{Nonlinear system theory: Another look at dependence}.
\newblock \bibinfo{journal}{Proc. Natl. Acad. Sci.} \bibinfo{volume}{102},
  \bibinfo{pages}{14150--14154}.
\bibitem[{Wu and Min(2005)}]{WuMin2005}
\bibinfo{author}{Wu, W.B.}, \bibinfo{author}{Min, M.}, \bibinfo{year}{2005}.
\newblock \bibinfo{title}{On linear processes with dependent innovations}.
\newblock \bibinfo{journal}{Stochastic Process. Appl.} \bibinfo{volume}{115},
  \bibinfo{pages}{939--958}.
\bibitem[{Zhang and Cheng(2018)}]{ZhangCheng2018}
\bibinfo{author}{Zhang, X.}, \bibinfo{author}{Cheng, G.}, \bibinfo{year}{2018}.
\newblock \bibinfo{title}{Gaussian approximation for high dimensional vector
  under physical dependence}.
\newblock \bibinfo{journal}{Bernoulli} \bibinfo{volume}{24},
  \bibinfo{pages}{2640--2675}.

\end{thebibliography}
\singlespacing\setstretch{1} \clearpage

\end{document}